\documentclass[12pt]{article}
\usepackage{amsmath,amssymb}
\newcommand{\beq}{\begin{eqnarray}}
\newcommand{\eeq}{\end{eqnarray}}
\newcommand{\rd}{\partial}

\newcommand{\bbC}{\mathbb{C}}
\newcommand{\bbP}{\mathbb{P}}

\newcommand{\calC}{\mathcal{C}}

\newcommand{\calU}{\mathcal{U}}
\newcommand{\calJac}{\mathcal{J}ac}
\newcommand{\calPic}{\mathcal{P}ic}
\newcommand{\Jac}{\mathrm{Jac}}
\newcommand{\Pic}{\mathrm{Pic}}
\newcommand{\Sym}{\mathrm{Sym}}
\newcommand{\SU}{\mathrm{SU}}
\begin{document}

\title{Hyperelliptic Integrable Systems on \\
K3 and Rational Surfaces}
\author{Kanehisa Takasaki\\
\normalsize Department of Fundamental Sciences, \\
\normalsize Faculty of Integrate Human Studies, Kyoto University\\
\normalsize Yoshida, Sakyo-ku, Kyoto 606-8501, Japan\\
\normalsize E-mail: takasaki@math.h.kyoto-u.ac.jp}
\date{}
\maketitle

\begin{abstract}
We show several examples of integrable systems 
related to special K3 and rational surfaces (e.g.,  
an elliptic K3 surface, a K3 surface given by a 
double covering of the projective plane, a rational 
elliptic surface, etc.).  The construction, based on 
Beauvilles's general idea, is considerably simplified 
by the fact that all examples are described by 
hyperelliptic curves and Jacobians.  This also enables 
to compare these integrable systems with more classical 
integrable systems, such as the Neumann system and the 
periodic Toda chain, which are also associated with 
rational surfaces.  A delicate difference between the 
cases of K3 and of rational surfaces is pointed out 
therein. 
\end{abstract}

\vfill

\begin{flushleft}
KUCP-0161\\
math.AG/0007073
\end{flushleft}

\newpage

\section{Introduction}

Some ten years ago, Beauville \cite{bib:Beauville} 
presented a construction of algebraically completely 
integrable Hamiltonian system (ACIHS) associated 
with a K3 surface.  Beauville's construction is 
based on Mukai's work \cite{bib:Mukai} on symplectic 
geometry of moduli spaces of sheaves on K3 surfaces.  
Beauville considered a special case of Mukai's 
moduli spaces, and discovered that this moduli space 
(with a symplectic structure in the sense of Mukai) 
is an ACIHS, namely, it is a symplectic variety with 
a Lagrangian fibration whose fibers are Abelian varieties.  
An interesting byproduct of his result is the fact 
that this ACIHS has another realization as a symmetric 
product of the K3 surface.   

Beauville's construction suggests a new approach 
to the issue of construction of integrable systems, 
namely, an approach from algebraic surfaces.  
Hurtubise \cite{bib:Hurtubise} indeed found 
such a general framework.  His result covers   
many classical integrable systems related to loop 
algebras \cite{bib:AvM,bib:AHH,bib:AHP,bib:RSTS}, 
Hitchin's integrable systems of Higgs pairs 
\cite{bib:Hitchin-higgs}, etc., as well as 
Beauville's integrable systems, on an equal 
footing.  Furthermore, Gorsky et al. 
\cite{bib:Go-Ne-Ru} pointed out a close 
relation between this approach and Sklyanin's 
separation of variables \cite{bib:Sklyanin}.  

As for the case of K3 surfaces, however, 
very few examples of Beauville's integrable 
systems seem to be {\it explicitly} known. 
Almost the only example in the literature is 
the case of quartics in $\bbP^3$ that Beauville 
mentioned in his paper to illustrate his 
construction.  (A more detailed description 
of this example can be found in the paper of 
Gorsky et al. \cite{bib:Go-Ne-Ru}) 

In this paper, we report a few examples of 
integrable systems related to Beauville's 
construction.  Two of them are direct application 
of Beauville's construction to an elliptic K3 
surface and a K3 surface given by a ramified 
double covering of $\bbP^2$.  Another example 
is an integrable system associated with 
a rational elliptic surface, which is also 
interesting in comparison with elliptic K3 
surfaces.  A common characteristic of these examples 
is that they are described by hyperelliptic curves 
and Jacobi varieties.  This enables us to compare 
these examples with more classical integrable 
systems that are related to hyperelliptic curves.

\section{Beauville's construction}

Following Beauville, we assume that the K3 
surface $X$ in consideration has an projective 
embedding or finite morphism $\phi: X \to \bbP^g$ 
by which a $g$-dimensional family of genus $g$ 
curves $C_u = X \cap \phi^{-1}(H_u)$ are cut out 
by hyperplanes $H_u$, $u \in (\bbP^g)^*$.  
In order to avoid delicate problems, we restrict 
the hyperplane parameters $u$ to a suitable open 
subset $\calU$ of $(\bbP^g)^*$.  We thus have 
a family $\calC \to \calU$ of curves fibered over 
$\calU$.  

The phase space of Beauville's integrable system 
is the associated relative Picard scheme 
$h: \calPic^g \to \calU$ of the degree $g$ 
component $\Pic^g(C_u)$ of the Picard group of 
$C_u$.  This relative Picard scheme, upon suitably 
compactified, can be identified with a moduli 
space of simple sheaves on $X$, so that, by 
Mukai's work, it has a symplectic structure.  
As Beauville discovered, $h$ is a Lagrangian 
fibration with respect to Mukai's symplectic 
structure.  Since each fiber $h^{-1}(u) = 
\Pic^g(C_u)$ is obviously an Abelian variety, 
the total space $\calPic^g$ of the relative 
Picard scheme becomes an ACIHS (see the review 
of Donagi and Markman \cite{bib:Do-Ma} for the 
notion of ACIHS).  

Another realization of this integrable system 
is given by the $g$-fold symmetric product $X^{(g)} 
= \Sym^g(X)$ equipped with a symplectic structure 
induced by the complex symplectic structure of $X$.  
In order to see the relation with the relative 
Picard scheme, recall that $\Pic^g(C_u)$ consists 
of the linear equivalence $[D]$ of an effective 
divisor $D = p_1 + \cdots + p_g$ of degree $g$.  
The unordered $g$-tuple $(p_1,\ldots,p_g)$ 
of points then becomes a point of $X^{(g)}$.  
Conversely, the $g$-tuple $(p_1,\ldots,p_g)$ 
of points in a {\it general position} uniquely 
determines a hyperplane $H_u \subset \bbP^g$ 
(hence a curve $C_u = X \cap H_u$) that passes 
through the $g$ points.  One can thus define 
the map 
\beq
\begin{array}{rcl}
    X^{(g)}_* & \longrightarrow & \calU \\
    (p_1,\ldots,p_g) & \longmapsto & u 
\end{array}
\eeq
from an open subset $X^{(g)}_*$ to $\calU$.  
This map can be lifted up to an open embedding 
\beq
\begin{array}{rcl}
    X^{(g)}_* & \longrightarrow & \calPic^g \\
    (p_1,\ldots,p_g) & \longmapsto & 
    (u,[p_1 + \cdots + p_g]) 
\end{array}
\eeq
by adding the linear equivalence class 
$[p_1 + \cdots + p_g] \in \Pic^g(C_u)$ as 
the second data.  Beauville observed that 
this is a symplectic mapping.  Thus, 
in particular, $X^{(g)}$ is fibered by 
the $g$-fold symmetric product $C_u^{(g)} 
= \Sym^g(C_u)$ of the curves $C_u$, and 
these fibers are Lagrangian.

\section{Hyperelliptic integrable system on 
elliptic K3 surface} 

We now apply Beauville's construction to 
an elliptic K3 surface $X \to \bbP^1$.  
The Weierstrass model of such an algebraic 
surface can be written, in affine coordinates, 
as 
\beq
    y^2 = z^3 + f(x) z + g(x), 
\eeq
where $x$ is an affine coordinate of $\bbP^1$ 
and $f(x)$ and $g(x)$ are polynomials of 
degree $8$ and $12$, respectively.  The 2-form 
\beq
    \omega = \frac{dz \wedge dx}{y}, 
\eeq
is nowhere-vanishing and holomorphic, so that 
$X$ becomes a complex symplectic surface.

\subsection{Curves, divisors and symmetric product}

We define a five-parameter family of curves $C_u$, 
$u = (u_1,\ldots,u_5)$, cut out from $X$ by the equation 
\beq
    z = P(x) = \sum_{k=1}^5 u_k x^{5-k}. 
\eeq
In other words, the curves $C_u$ are hyperelliptic 
curves of genus $5$ defined by the equation 
\beq
    y^2 = P(x)^3 + f(x) P(x) + g(x). 
\eeq
Note that the genus and the number of parameters 
of curves coincide.  Furthermore, these curves 
are hyperplane sections of the following 
finite morphism $\phi: X \to \bbP^5$: 
\beq
\begin{array}{rcl}
    X & \longrightarrow & \bbP^5 \\
    (x,y,z) & \longmapsto & (-z:x^4:x^3:x^2:x^1:x:1)
\end{array}
\eeq
Although we have used the affine coordinates for 
simplicity, it is not difficult to projectivize 
this definition.  The five parameters $u_1,\dots,u_5$ 
thus turn out to be affine coordinates of $\bbP^5$. 

The mapping $X^{(5)}_* \to \calU$ takes a particularly 
simple form in the present setting.  Given an unordered 
$5$-tuple $(p_1,\ldots,p_5)$ of points of $X$, 
let $(x_j,y_j,z_j)$ denote the affine coordinates 
of the $j$-th point $p_j$.  The curve $C_u$ passes 
through the five points if and only if the equations 
\beq
    P(x_j) = z_j 
\eeq
are satisfied for $j = 1,\ldots,5$.  The polynomial 
$P(x)$ is uniquely determined by these conditions if 
the $x_j$'s are distinct: $x_j \not= x_k$  
($j \not= k$).  In fact, the Lagrange interpolation 
formula provides an explicit expression of $P(x)$: 
\beq
    P(x) = \sum_{j=1}^5 z_j 
      \prod_{k\not=j} \frac{x - x_k}{x_j - x_k}. 
\eeq
This determines $u_j$'s as rational functions 
of $(x_1,\ldots,x_5,z_1,\ldots,z_5)$, e.g., 
\beq
    u_1 = \sum_{j=1}^5 z_j / \prod_{k\not=j}(x_j - x_k). 
\eeq
These functions are to be a commuting set of 
Hamiltonians.

\subsection{Abel-Jacobi mapping and symplectic form}

The symplectic structure of $X^{(5)}$ is 
defined by the 2-form 
\beq
    \Omega = 
    \sum_{j=1}^5 \frac{dz_j \wedge dx_j}{y_j}, 
\eeq
where $(x_j,y_j,z_j)$ denote the affine 
coordinates of the $p_j$'s.  According to 
Beauville, this symplectic structure is 
mapped to Mukai's symplectic structure on 
$\calPic^5$.  

Following Hurtubise, we now slightly modify 
Beauville's construction:  We shift the divisors 
$p_1 + \cdots + p_5$ to $p_1 + \cdots + p_5 
- 5p_\infty$, where $p_\infty$  is the point 
$(x,y,z) = (\infty,\infty,\infty)$ of $C_u$, 
and consider the relative Jacobian $\calJac 
\to \calU$ of the Jacobi varieties $\Jac(C_u)$ 
fibered over $\calU$ rather than the relative Picard.  
The Abel-Jacobi mapping on $C_u$ then induces 
a mapping $X^{(5)}_* \to \calJac$, which plays 
the role of Beauville's mapping.  

The Abel-Jacobi mapping sends the symmetric 
product $C_u^{(5)}$ onto the Jacobi variety 
$\Jac(C_u) = \bbC^5/L_u$ realized as a complex 
torus with period lattice $L_u$: 
\beq
\begin{array}{rcl}
    C_u^{(5)} & \longrightarrow & 
    \Jac(C_u) = \bbC^5/L_u \\
    (p_1,\ldots,p_5) & \longmapsto & 
    (\psi_1,\ldots,\psi_5) 
\end{array}
\eeq
In order to define this mapping, one has to 
choose a basis of holomorphic differentials 
on $C_u$.  We take the standard holomorphic 
differentials 
\beq
    \sigma_k = \frac{x^{5-k}dx}{y} 
    \quad (k = 1,\ldots,5).  
\eeq
The $\psi_k$'s are then given by a sum of 
integrals of the $\sigma_k$'s: 
\beq
    \psi_k = \sum_{j=1}^5 
      \int_{(\infty,\infty)}^{(x_j,y_j)} 
      \frac{x^{5-k}dx}{y}
\eeq

Now $(u_1,\ldots,u_5,\psi_1,\ldots,\psi_5)$ give 
a local coordinate system on the relative Jacobian 
$\calJac$.  Remarkably, $\Omega$ turns out to take 
a canonical form in these coordinates: 
\beq\label{eq:Omega-dudpsi}
    \Omega = \sum_{k=1}^5 du_k \wedge d\psi_k. 
\eeq
In other words, the ``toroidal'' coordinates $\psi_k$' 
are nothing but the canonical conjugate variables 
of the hyperplane parameters $u_k$.

\subsection{Proof of (\ref{eq:Omega-dudpsi})}

Since 
\[
    z_j = P(x_j) = \sum_{k=1}^5 u_k x_j^{5-k}, 
\]
$\Omega$ can be written 
\beq
    \Omega = \sum_{j,k=1}^5 
      \frac{du_k \wedge x_j^{5-k}dx_j}{y_j}. 
\eeq
It is convenient choose 
$(x_1,\ldots,x_5,\psi_1,\ldots,\psi_5)$, 
rather than 
$(u_1,\ldots,u_5,\psi_1,\ldots,\psi_5)$, 
as independent variables here.  The total 
differential of $d\psi_k$ can be now expressed as 
\[
    d\psi_k = 
      \sum_{j=1}^5 \frac{\rd\psi_k}{\rd x_j}dx_j 
    + \sum_{\ell=1}^5 \frac{\rd\psi_k}{\rd u_\ell}du_\ell, 
\]
so that 
\beq
    \sum_{k=1}^5 du_k \wedge d\psi_k 
    = \sum_{j,k=1}^5 
        \frac{\rd\psi_k}{\rd x_j}du_k \wedge dx_j 
    + \sum_{k,\ell=1}^5 
         \frac{\rd\psi_k}{\rd u_\ell}du_k \wedge du_\ell. 
\eeq
Differentiating the explicit definition 
\[
    \psi_k = \sum_{j=1}^5 
      \int_{\infty}^{x_j}
      \frac{x^{5-k}dx}
      {\Bigl(P(x)^3 + f(x)P(x) + g(x)\Bigr)^{1/2}} 
\]
of the Abel-Jacobi mapping gives 
\beq
    \frac{\rd\psi_k}{\rd x_j} 
    &=& \frac{x_j^{5-k}}{y_j}, 
    \\
    \frac{\rd\psi_k}{\rd u_\ell}
    &=& \sum_{j=1}^5 \int_\infty^{x_j} 
        \frac{x^{5-k}\cdot x^{5-\ell} 
              \Bigl(3P(x)^2 + f(x)\Bigr)dx} 
        {2\Bigl(P(x)^3 + f(x)P(x) + g(x)\Bigr)^{3/2}}. 
\eeq
In particular, 
\beq\label{eq:dpsidu} 
    \frac{\rd\psi_k}{\rd u_\ell} 
  = \frac{\rd\psi_\ell}{\rd u_k}, 
\eeq
therefore 
\beq
    \sum_{k,\ell=1}^5 
    \frac{\rd\psi_k}{\rd u_\ell} du_k \wedge du_\ell 
    = 0. 
\eeq
This implies the equality 
\beq
    \sum_{k=1}^5 du_k \wedge d\psi_k 
  = \sum_{j,k=1}^5 \frac{x_j^{5-k}}{y_j} du_k \wedge dx_j 
  = \Omega, 
\eeq
thus completing the proof of (\ref{eq:Omega-dudpsi}).

\subsection{Remarks}

\paragraph*{1.} 
Essentially the same expression of 
the symplectic form is stated in 
Hurtubise's paper \cite{bib:Hurtubise} 
in a more general form.  This expression 
is thus not specific to the case of 
hyperelliptic curves, but rather universal.  

\paragraph*{2.} 
(\ref{eq:Omega-dudpsi}) is a concise expression 
of the fact that the relative Jacobian $\calJac 
\to \calU$ is an ACIHS.  Firstly, the involutivity 
(Poisson-commutativity) of the ``Hamiltonians'' $u_k$ 
is obvious from this expression of $\Omega$.  
Accordingly, the intersection of the level surfaces 
of the $u_k$'s, which are nothing but the Jacobi 
varieties $\Jac(C_u)$, are Lagrangian subvarieties.  
Lastly, these Hamiltonians generate a commuting set 
of {\it linear} flows on these complex tori. 
The Hamiltonian flow of $u_j$ is indeed given 
by 
\beq
    \psi_k(t) = \psi_k(0) + \delta_{jk}t. 
\eeq

\paragraph*{3.} 
$\Omega$ should be invariant under the shift 
\beq
    \psi_k \ \longrightarrow \ 
    \psi_k + e_k 
\eeq
by any element $e = (e_1,\ldots,e_5)$ of 
the period lattice $L_u$.  This implies 
the nontrivial equations 
\beq
    \sum_{k=1}^5 du_k \wedge de_k = 0 
\eeq
among the components of the vector $e$ 
(which depends on $u$).  This is an avatar 
of the ``cubic condition'' of Donagi 
and Markman \cite{bib:Do-Ma}.  Actually, 
we can verify these conditions directly: 
Note that $e_k$ is given by the period integral 
\beq
    e_k = \int_\gamma \frac{x^{5-k}dx}{y} 
\eeq
along a cycle $\gamma$.  From this 
expression, one can derive the equality 
\beq
    \frac{\rd e_k}{\rd u_\ell} 
  = \frac{\rd e_\ell}{\rd u_k} 
\eeq
in the same way as the calculations in the 
proof of (\ref{eq:Omega-dudpsi}).  

\paragraph*{4.} 
A one-parameter family of hyperelliptic curves 
of genus $5$ on the elliptic K3 surface appears 
in the work of Fock et al. on duality of integrable 
systems \cite{bib:Fo-Go-Ne-Ru}.  In fact, their 
hyperelliptic curves are the subfamily of ours 
with $u_1 = u_2 = u_3 = u_4 = 0$.

\section{Hyperelliptic integrable system on 
non-elliptic K3 surface}

We here present an example associated with a 
double covering of $\bbP^2$.   Let $(x,z)$ be 
affine coordinates of $\bbP^2$ and consider an 
algebraic surface $X$ defined by the equation 
\beq
    y^2 = f(x,z), 
\eeq
where $f(x,z)$ is a polynomial of degree six 
whose zero-locus $f(x,z) = 0$ is a nonsingular 
sextic curve on $\bbP^2$.  This is a K3 surface 
that appears in Mukai's work as an interesting 
example \cite{bib:Mukai}.   A nowhere-vanishing 
holomorphic 2-form is given by 
\beq
    \omega = \frac{dz \wedge dx}{y}. 
\eeq

A two-parameter family of curves $C_u$, 
$u = (u_1,u_2)$, can be cut out from $X$ 
by the equation 
\beq
    z = P(x) = u_1 x + u_2. 
\eeq
In other words, $C_u$ is the hyperelliptic 
curve of genus $2$ defined by the equation 
\beq
    y^2 = f(x, u_1 x + u_2). 
\eeq
Note that these curves are nothing but the 
hyperplane sections associated with the 
double covering $\phi: X \to \bbP^2$.  

The rest of the construction is fully parallel 
to the case of the elliptic K3 surface.  
The integrable system is realized on (an open 
subset of) the $2$-fold symmetric product 
$X^{(2)}$ or on the relative Jacobian $\calJac$.  
Let $(p_1,p_2)$ denote an unordered pair of 
points of $X$, $(x_j,y_j,z_j)$ the coordinates 
of $p_j$, and $\Omega$ the symplectic form 
induced on $X^{(2)}$: 
\beq
    \Omega = \sum_{j=1,2} \frac{dz_j \wedge dx_j}{y_j}. 
\eeq
The Abel-Jacobi mapping 
\beq
\begin{array}{rcl}
  C_u^{(2)} & \longrightarrow & \Jac(C_u) = \bbC^2/L_u \\
  (p_1,p_2) & \longmapsto & (\psi_1,\psi_2) 
\end{array}
\eeq
is now given by 
\beq
   \psi_1 = \sum_{j=1,2} 
     \int_{(\infty,\infty)}^{(x_j,y_j)} \frac{xdx}{y}, 
   \quad 
   \psi_2 = \sum_{j=1,2} 
     \int_{(\infty,\infty)}^{(x_j,y_j)} \frac{dx}{y}. 
\eeq
The symplectic form $\Omega$ takes the canonical form 
\beq
    \Omega = \sum_{k=1,2} du_k \wedge d\psi_k 
\eeq
in the coordinates $(u_1,u_2,\psi_1,\psi_2)$.

\section{Hyperelliptic integrable system on 
rational surface}

\subsection{Integrable system associated with 
rational elliptic surface}

If the polynomials $f(x)$ and $g(x)$ in the 
Weierstrass model of elliptic K3 surfaces 
are replaced by polynomials of degree $4$ and 
$6$, resplectively, the outcome is an rational 
elliptic surface.  The 2-form 
\beq
    \omega = \frac{dz \wedge dx}{y}
\eeq
is nowhere-vanishing and holomorphic in the 
{\it affine part} of the surface, having poles 
at the compactification divisor at infinity.  
Let us examine this case in comparison with 
the case of K3 surfaces.  

Simple power counting suggests that the most 
natural choice of $P(x)$ in this case will be 
a polynomial of degree $2$: 
\beq
    P(x) = u_0 x^2 + u_1 x + u_2. 
\eeq
The equation $z = P(x)$ then defines a hyperelliptic 
curve of genus $2$ in $X$: 
\beq
    y^2 = P(x)^3 + f(x) P(x) + g(x). 
\eeq
Ganor \cite{bib:Ganor} considered the same family 
of curves in the context of string theory.  
The foregoing construction, however, does not work 
literally, firstly because the number of parameters 
(three) and the genus of the curves (two) do not match.  

A correct prescription is to fix $u_0$ to 
a constant $c$, and to take the two-parameter 
family of curves $C_u$, $u = (u_1,u_2)$, cut out 
from $X$ by the equation 
\beq
    z = P(x) = c x^2 + u_1 x + u_2. 
\eeq
In other words, $u_0 = c$ should be treated as 
a {\it Casimir function} on a Poisson manifold. 

Now we consider $X^{(2)}$ (rather than $X^{(3)}$ 
that Ganor argued) with the symplectic form 
\beq
    \Omega = \sum_{k=1,2} \frac{dz_j \wedge dx_j}{y_j}, 
\eeq
where $(x_j,y_j,z_j)$ are the coordinates of an 
unordered pair $(p_1,p_2)$ of points on $X$.  
The Abel-Jacobi mapping $(p_1,p_2) \mapsto 
(\psi_1,\psi_2)$  is given by 
\beq
    \psi_1 = \sum_{j=1,2} 
      \int_{(\infty,\infty)}^{(x_j,y_j)} \frac{xdx}{y}, 
    \quad 
    \psi_2 = \sum_{j=1,2} 
      \int_{(\infty,\infty)}^{(x_j,y_j)} \frac{dx}{y}, 
\eeq
and sends $\Omega$ the canonical form: 
\beq
    \Omega = \sum_{k=1,2} du_k \wedge d\psi_k.  
\eeq

Thus the situation for the rational elliptic surfaces 
is slightly different from the case of K3 surfaces.  
Namely, we have to fix the leading coefficient of 
$P(x)$ to a constant $c$ so as to match the number 
of parameters and the genus of curves.

\subsection{Neumann system and rational surface}

It is instructive to compare the foregoing 
construction with the Neumann system.  
The Neumann system is a classical integrable 
system that was solved in the 19th century 
\cite{bib:Neumann} using the theory of 
hyperelliptic integrals developed in those 
days, and revived by Moser \cite{bib:Moser} 
from a modern point of view. 

The hyperelliptic curves of the Neumann system 
take the special form 
\beq
    y^2 = P(x) Q(x),  
\eeq
where $P(x)$ and $Q(x)$ are polynomials of 
degree $N$ and $N+1$ of the form 
\beq
    P(x) = x^N + \sum_{k=1}^N u_k x^{N-k}, 
    \quad 
    Q(x) = \prod_{n=1}^{N+1} (x - c_n), 
\eeq
and $N$, the genus of the curves, is equal 
to the degrees of freedom of the Neumann 
system.  The coefficients $u_k$ of $P(x)$ 
are integrals of motion in involution.  
The $c_n$'s are non-dynamical structural 
constants (i.e., Casimir functions) of the 
system.  

One will soon notice that a rational surface 
is hidden behind these hyperelliptic curves. 
This rational surface $X$ is defined by the 
equation  
\beq
    y^2 = z Q(x). 
\eeq
The 2-form 
\beq
    \omega = \frac{dz \wedge dx}{y} 
\eeq
is nowhere-vanishing and holomorphic in 
the affine part of $X$.  The curves $C_u$ 
cut out from $X$ by the equation $z = P(x)$ 
are exactly the aforementioned hyperelliptic 
curves.  It is not difficult to confirm that 
the construction based on the symmetric product 
$X^{(N)}$ or the relative Jacobian $\calJac$ 
does reproduce the classical method for solving 
the Neumann system.  

It is interesting to note that the leading 
coefficient of $P(x)$, which is now set to $1$, 
is also a kind of structural constant.  This 
constant is related to the radius of an 
$N$-dimensional sphere on which the Neumann 
system is realized as a mechanical system.  
If the radius takes a different value, 
the leading coefficient of $P(x)$ accordingly 
varies as 
\[
    P(x) = cx^N + \sum_{k=1}^N u_k x^{N-k}. 
\]
This is the same situation that we have 
encountered in the case of rational elliptic 
surfaces.

\section{Concluding Remarks} 

We have shown several integrable systems 
that illustrate Beauville's construction.  
The construction is considerably simplified 
by the fact that the curves arising therein 
are hyperelliptic.   Nevertheless, the most 
essential part, i.e., the two expressions 
of $\Omega$ on the symmetric product $X^{(g)}$ 
and on the relative Jacobian $\calJac$, are 
rather universal and can be extended to a 
non-hyperelliptic case.  Let us stress that 
the coordinates $(x_j,y_j,z_j)$ of the unordered 
$g$-tuple $(p_1,\ldots,p_g) \in X^{(g)}$  
correspond to ``separated variables'' in 
Sklyanin's separation of variables. 

We have pointed out a delicate difference 
between the cases of K3 and of rational surfaces.  
Namely, some of the hyperplane parameters $u_k$ 
in the latter case have to be interpreted as 
a Casimir rather than an integral of motion.  
This issue deserves to be pursued further.  

The present approach from algebraic surfaces can 
be a promising alternative to the conventional 
method based on a Lax representation.  We have 
seen such a possibility in the case of the 
Neumann system.  Let us mention another example, 
which is related to Witten's eleven-dimensional 
interpretation  \cite{bib:Witten} of the 
Seiberg-Witten curves for four-dimensional 
$N=2$ supersymmetric QCD. 

The Seiberg-Witten curves for $\SU(N_c)$ 
gauge theory coupled to $N_f$ ($N_f \le 2N_c$) 
hypermultiplets are hyperelliptic curves of 
genus $g = N_c - 1$ defined by the equation 
\beq
    y^2 - P(x)y + Q(x) = 0. 
\eeq
Here $P(x)$ is a polynomial of the form 
\beq
    P(x) = x^{N_c} + \sum_{k=2}^{N_c} u_k x^{N_c-k},  
\eeq
and $Q(x)$ a polynomial of degree $N_f$ of 
the form 
\beq
    Q(x) = \Lambda^{2N_c - N_f} 
           \prod_{\ell=1}^{N_f} (x + m_\ell). 
\eeq
The $u_k$'s are the ``moduli'' of the 
Seiberg-Witten curves, which we interpret as 
(part of) the parameters of hyperplane 
sections.  If $N_f = 0$, $Q(x)$ reduces to 
a constant, and the curve is exactly the 
the spectral curve of the $N_c$-periodic 
Toda chain.  

In Witten's interpretation, these curves 
are embedded into the {\it affine} rational 
surface $X$ defined by the equation 
\beq
    yz = Q(x), 
\eeq
A nowhere-vanishing holomorphic 2-form 
on $X$ is given by 
\beq
    \omega = \frac{dy \wedge dx}{y} 
           = \frac{dx \wedge dz}{z}. 
\eeq
In fact, this complex surface is a realization 
of the multi-Taub-NUT space (gravitational 
instanton) in one of its $2$-parameter family 
of complex structures \cite{bib:Hitchin-graviton}.  
As noted by Nakatsu et al. \cite{bib:Nakatsu-etal} 
and de Boer et al. \cite{bib:deBoer-etal}, 
the Seiberg-Witten curves are cut out from $X$ 
by the equation 
\beq
    y + z = P(x). 
\eeq

One can now repeat the foregoing construction with 
a slightest modification to obtain an integrable 
system realized on $X^{(g)}$ or on $\calJac$. 
Presumably, this integrable system  will be 
equivalent to the spin chain that Gorsky et al. 
\cite{bib:Go-Ma-Mi-Mo} proposed as an 
integrable system for these Seiberg-Witten  
curves.  Our construction, however, is more 
direct and seemingly more natural in the context 
of brane theory \cite{bib:Go-Ne-Ru,bib:Witten}.

\subsection*{Acknowledgements}
I am grateful to Yasuhiko Yamada for 
useful comments.  This work is partly supported 
by the Grant-in-Aid for Scientific Research 
(No. 12640169), the Ministry of Education, 
Science and Culture.


\end{document}